\documentclass[doublespacing]{article}

\RequirePackage[OT1]{fontenc}

\RequirePackage[amsthm,amsmath,natbib]{imsart}
\RequirePackage{imsart}
\usepackage{amsmath,amssymb}
\startlocaldefs
\endlocaldefs

\begin{document}
%Synmboles grecs
\newcommand{\dd}{\delta}
\newcommand{\e}{\epsilon}
\newcommand{\alp}{\alpha}
\newcommand{\lab}{\lambda}
\newcommand{\gam}{\gamma}
\newcommand{\sig}{\sigma}
\newcommand{\tht}{\theta}
\newcommand{\Lab}{\Lambda}
\newcommand{\Gam}{\Gamma}
\newcommand{\Sig}{\Sigma}
\newcommand{\T}{\Theta}

%Sommes

\newcommand{\sli}{\sum\limits}
\newcommand{\sliin}{\sum\limits_{i=1}^n}
\newcommand{\sliid}{\sum\limits_{i=1}^d}
\newcommand{\sliik}{\sum\limits_{i=1}^k}
\newcommand{\sliiN}{\sum\limits_{i=1}^N}
\newcommand{\slijk}{\sum\limits_{j=1}^k}
\newcommand{\proliin}{\prod\limits_{i=1}^n}
\newcommand{\proliid}{\prod\limits_{i=1}^d}
\newcommand{\proliik}{\prod\limits_{i=1}^k}
\newcommand{\proliiN}{\prod\limits_{i=1}^N}
\newcommand{\prolijk}{\prod\limits_{j=1}^k}
\newcommand{\ili}{\int\limits}
\newcommand{\proli}{\prod\limits}
\newcommand{\bculi}{\bigcup\limits}
\newcommand{\bcali}{\bigcap\limits}

%Convergences

\newcommand{\ls}{\limsup}
\newcommand{\li}{\liminf}
\newcommand{\limn}{\lim_{n\rightarrow\infty}\;}
\newcommand{\limk}{\lim_{k\rightarrow\infty}\;}
\newcommand{\lsn}{\limsup_{n\rightarrow\infty}}
\newcommand{\lin}{\liminf_{n\rightarrow\infty}}
\newcommand{\lsk}{\limsup_{k\rightarrow\infty}}
\newcommand{\lik}{\liminf_{k\rightarrow\infty}}
\newcommand{\rar}{\rightarrow}
\newcommand{\lar}{\leftarrow}
\newcommand{\cvps}{\rightarrow_{p.s.}\;}
\newcommand{\cvpr}{\rightarrow_{P}\;}
\newcommand{\cvloi}{\rightarrow_{\mcal{L}}\;}
\newcommand{\kif}{k\rightarrow\infty}
\newcommand{\nif}{n\rightarrow\infty}

%Parenthèses

\newcommand{\aoo}{\Big\{}
\newcommand{\aff}{\Big\}}
\newcommand{\coo}{\Big [}
\newcommand{\cff}{\Big]}
\newcommand{\poo}{\Big (}
\newcommand{\pff}{\Big)}
\newcommand{\po}{\big (}
\newcommand{\pf}{\big)}
\newcommand{\ao}{\big \{}
\newcommand{\af}{\big \}}
\newcommand{\co}{\big [}
\newcommand{\cf}{\big ]}
\newcommand{\pooo}{\bigg (}
\newcommand{\pfff}{\bigg)}
\newcommand{\aooo}{\bigg \{}
\newcommand{\afff}{\bigg \}}
\newcommand{\cooo}{\bigg [}
\newcommand{\cfff}{\bigg ]}
\newcommand{\poooo}{\Bigg (}
\newcommand{\pffff}{\Bigg)}
\newcommand{\aoooo}{\Bigg \{}
\newcommand{\affff}{\Bigg \}}
\newcommand{\coooo}{\Bigg [}
\newcommand{\cffff}{\Bigg ]}
\newcommand{\mmi}{\mid\mid}
\newcommand{\mmmi}{\mid\mid\mid}
\newcommand{\mmii}{\mid\mid_\infty}
\newcommand{\gMid}{\Bigg |}
\newcommand{\Mid}{\Big |}
\newcommand{\Mmi}{\Big|\Big|}
\newcommand{\Mmii}{\Big|\Big|_\infty}

%Tableaux triangulaires, suites

\newcommand{\eni}{\po \e_{n,i}\pf_{n\geq1, i \le p_n}}
\newcommand{\tab}{_{n\geq1,\;i\le p_n}}
\newcommand{\atc}{n\geq 1,\;i\le p_n}
\newcommand{\suite}{_{n\geq 1}}

%Sous suites

\newcommand{\nk}{n_k}
\newcommand{\nkm}{n_{k-1}}
\newcommand{\nkk}{n_{k+1}}
\newcommand{\nkd}{n_k}
\newcommand{\nkkd}{n_{k+1}}

%Paramètres LLI

\newcommand{\hnk}{h_{n_k}}
\newcommand{\bnk}{b_{n_k}}
\newcommand{\bnkk}{b_{n_{k+1}}}
\newcommand{\hkl}{{h_{n_k,l}}}
\newcommand{\hnl}{{h_{n,l}}}
%Mathbb

\newcommand{\EEE}{\mathbb{E}}
\newcommand{\NNN}{\mathbb{N}}
\newcommand{\PPP}{ \mathbb{P}}
\newcommand{\CCC}{\mathbb{C}}
\newcommand{\KKK}{\mathbb{K}}
\newcommand{\RRR}{\mathbb{R}}
%Mathcal

\newcommand{\FF}{\mathcal{F}}
\newcommand{\TT}{\mathcal{T}}
\newcommand{\GG}{\mathcal{G}}
\newcommand{\BBGG}{\mathcal{B}(\mathcal{G})}
\newcommand{\CC}{\mathcal{C}}
\newcommand{\KK}{\mathcal{K}}
\newcommand{\SSS}{\mathcal{S}}
\newcommand{\BB}{\mathcal{B}}
\newcommand{\PP}{\mathcal{P}}
\newcommand{\HH}{\mathcal{H}}
\newcommand{\NN}{\mathcal{N}}
\newcommand{\MM}{\mathcal{M}}
\newcommand{\DD}{\mathcal{D}}

%Surlignage
\newcommand{\ovg}{\overline{g}}
\newcommand{\ovh}{\overline{h}}
\newcommand{\ovE}{\overline{E}}
\newcommand{\ovH}{\overline{H}}
\newcommand{\ovI}{\overline{I}}
\newcommand{\ovJ}{\overline{J}}
\newcommand{\ovK}{\overline{K}}
\newcommand{\ovR}{\overline{R}}
\newcommand{\ovFF}{\overline{\mathcal{F}}}
\newcommand{\ovPPP}{\overline{\mathbb{P}}}

%Tildage
\newcommand{\wt}{\widetilde}
\newcommand{\cc}{\widetilde{c}}
\newcommand{\wtG}{\widetilde{G}}
\newcommand{\wtf}{\widetilde{f}}
\newcommand{\wth}{\mathfrak{h}}
\newcommand{\wtn}{\widetilde{n}}
\newcommand{\wtv}{\widetilde{v}}
\newcommand{\wtA}{\widetilde{A}}
\newcommand{\wtC}{\widetilde{C}}
\newcommand{\wtE}{\widetilde{E}}
\newcommand{\wtF}{\widetilde{F}}
\newcommand{\wtI}{\widetilde{I}}
\newcommand{\wtK}{\widetilde{K}}
\newcommand{\wtJ}{\widetilde{J}}
\newcommand{\wtN}{\widetilde{N}}
\newcommand{\wtP}{\widetilde{P}}
\newcommand{\wtFF}{\widetilde{\mathcal{F}}}
\newcommand{\wtPPP}{\widetilde{\mathbb{P}}}

%Divers

\newcommand{\lb}{\newline}
\newcommand{\inv}{\frac{1}}
\newcommand{\indep}{{\bot}\kern-0.9em{\bot}}
\newcommand{\beq}{\begin{equation} }
\newcommand{\eeq}{\end{equation} }
\newcommand{\mcal}{\mathcal}
\newcommand{\sq}{\sqrt}
\newcommand{\srl}{\stackrel}
\newcommand{\wap}{(\Omega,\mathcal{A},\rm I\kern-2pt P)}
\newcommand{\vk}{\vskip10pt}
\newcommand{\norm}{\mid\mid \cdot \mid\mid}
\newcommand{\nono}{\nonumber}

\newcommand{\Cov}{\mathrm{Cov}}
\newcommand{\Var}{\mathrm{Var}}
\newcommand{\FFGG}{\FF\times\GG}

\newtheorem{theo}{Theorem}
\newtheorem{fact}{Fact}
\newtheorem{ptheo}{Preuve du théorème}[section]
\newtheorem{lem}{Lemma}[section]
\newtheorem{plem}{Proof of Lemma}[section]
\newtheorem{prop}{Propriété}[section]
\newtheorem{preuveprop}{Preuve de la propriété}[section]
\newtheorem{defi}{Definition}[section]
\newtheorem{propo}{Proposition}[section]
\newtheorem{popo}{Proof of Proposition}[section]
\newtheorem{coro}{Corollary}[theo]
\newtheorem{pcoro}{Preuve du corollaire}[theo]
\newtheorem{ineg}{Inégalité}[section]
\newtheorem{pineg}{Preuve de l'inégalité}[section]
\newtheorem{rem}{Remark}[section]
\numberwithin{equation}{section}
\begin{frontmatter}

% "Title of the paper"
\title{Some uniform in bandwidth functional results for the tail uniform empirical and
quantile processes} \runtitle{Uniform in bandwidth results for
tail empiricals}

% indicate corresponding author with \corref{}
%\author{\fnms{John} \snm{Smith}\corref{}\ead[label=e1]{smith@foo.com}\thanksref{t1}}
% \thankstext{t1}{Thanks to somebody}
% \address{line 1\\ line 2\\ printead{e1}}
 \affiliation{Institute of Statistics, Catholic University of Louvain}

\author{\fnms{Davit} \snm{Varron}\ead[label=e1]{varron@stat.ucl.ac.be}}
\address{\printead{e1}}
%\and
%\author{\fnms{???} \snm{???}\ead[label=e2]{???}}
%\address{\printead{e2}}
%\affiliation{???}

\runauthor{D. Varron}

\begin{abstract}
For fixed $t\in [0,1)$ and $h>0$, consider the local uniform
empirical process
$$\DD_{n,h,t}(s):=n^{-1/2}\coo\sliin 1_{[t,t+hs]}(U_i)-hs\cff,\;s\in
[0,1],$$ where the $U_i$ are independent and uniformly distributed
on $[0,1]$. We investigate the functional limit behaviour of
$\DD_{n,h,t}$ uniformly in $\wth_n\le h\le h_n$ when
$n\wth_n/\log\log n\rar \infty$ and $h_n\rar 0$.
\end{abstract}



\end{frontmatter}
\section{Introduction}
Let $(U_i)_{i\geq 1}$ be an independent, identically distributed
(i.i.d.) sequence of random variables that are uniformly
distributed on $[0,1]$. Define the empirical distribution function
based on $(U_1,\ldots,U_n)$ by $F_n(t):=n^{-1}\sharp\{1\le i\le
n,\;U_i\le t\},\;t\in [0,1]$ and denote by $F_n^{\leftarrow}(t)$
the left-continuous inverse of $F_n$, namely
$F_n^{\leftarrow}(t):=\inf\{s\geq 0,\;F_n(s)\geq t\}$. We also
define the empirical (resp. quantile) process by
$\alp_n(t):=\sqrt{n}(F_n(t)-t),t\in[0,1]$ (resp.
$\beta_n(t):=\sqrt{n}(F_n^{\leftarrow}(t)-t),\;t\in[0,1]$). The
framework of this paper is the almost sure behaviour of the local
empirical and quantile processes. Namely, given $t\in[0,1)$ we
focus on studying the following processes, as $\nif$ and $h\rar
0$.
\begin{align}
\DD_{n,h,t}(s):=&\alp_n(t+h
s)-\alp_n(t),\;s\in[0,1],\label{Deltanht}\\
\DD'_{n,h,t}(s):=&\beta_n(t+h
s)-\beta_n(t),\;s\in[0,1].\label{Deltpanht}
\end{align}
Mason (1988) was the first to establish a functional law of the
iterated logarithm for the local empirical process (see also
Einmahl and Mason (1997) for a generalization of this result to
empirical processes indexed by functions). To cite this result, we
need to introduce some further notations first. Write $\log_2(u):=
\log(\log (u\vee 3)).$ We say that a sequence $(h_n)_{n\geq 1}$ of
strictly positive constants satisfies the local strong invariance
conditions when, ultimately as $\nif$, \beq h_n\downarrow
0,\;nh_n\uparrow \infty,\;nh_n/\log_2 n \rar
\infty.\label{CRS}\eeq Given a sequence $(x_n)_{n\geq 1}$ of
elements of a metric space $(E,d)$, we say that $x_n\leadsto K$
when $K$ is non void and coincides with the set of all cluster
points of $(x_n)\suite$. In our framework, $(E,d)$ is the space
$B([0,1])$ of all real bounded CADLAG trajectories on $[0,1]$,
endowed with the usual sup norm, namely $\mmi g \mmi:=\sup\{\mid
g(s)\mid,\;s\in [0,1]\}$. Consider the space $AC[0,1]$ of all
absolutely continuous functions on $[0,1]$. For any $g\in
AC[0,1]$, we define the usually called Hilbertian norm of $g$ as
\beq\mmi g\mmi_H^2:= \ili_0^1 {\dot g}^2(x)dx\label{normeH},\eeq
where $\dot g$ is any version of the derivative of $g$ with
respect to the Lebesgue measure. The usually called Strassen ball
can be defined as follows: \beq \SSS:=\aoo g\in AC([0,1]),\;
g(0)=0, \mmi g\mmi_H\le 1\aff.\label{SSS}\eeq As a corollary of a
strong approximation result, Mason (1988) showed that, given a
sequence $(h_n)\suite$ fulfilling (\ref{CRS}) and given $t\in
[0,1)$, we have, almost surely \beq \frac{\DD_{n,h_n,t}}{(2
h_n\log_2 n)^{1/2}}\leadsto \SSS\label{Masonloc}\eeq In the
particular case where $t=0$, Einmahl and Mason (1988) showed that
$\DD'_{n,h_n,t}$ also satisfies (\ref{Masonloc}). They showed that
result by making use of a local Bahadur Kiefer representation (see
their Theorem 5). The almost sure limit behavior of
$\DD'_{n,h_n,t}$ when $t\in (0,1)$ has been investigated by
Deheuvels (1997), who showed that the above mentioned process may
obey functional limit laws that are different from
(\ref{Masonloc}). The aim of the present paper is the following:
given two sequences $\wth_n<h_n$ fulfilling (\ref{CRS}), does
(\ref{Masonloc}) still hold uniformly in $\wth_n\le h\le h_n$?
Namely, do we have almost surely
\begin{align}
\limn \sup_{\wth_n\le h\le h_n}\;\inf_{g\in \SSS}\Mmi
\frac{\DD_{n,h,t}}{(2h\log_2 n)^{1/2}}- g\Mmi =&0,\label{fau}\\
\forall g \in \SSS,\;\lin \sup_{\wth_n\le h\le h_n}\Mmi
\frac{\DD_{n,h,t}}{(2h\log_2 n)^{1/2}}- g\Mmi =0\label{faux}\;?
\end{align}
The remainder of this paper is organised as follows. In \S
\ref{deux}, we state our main results on $\DD_{n,h,t}$. We then
show how this results lead to a local Bahadur-Kiefer type
representation that holds uniformly in $h$. The proofs of our main
results follow in \S\ref{trois}, \ref{quatre} and \ref{cinq}.
\section{Mains results}\label{deux}
Our first result is a weaker form of assertion (\ref{fau}).
\begin{theo}\label{T1}
Let $(h_n)\suite$ and $(\wth_n)\suite$ be two sequences satisfying
(\ref{CRS}) as well as $ \wth_n<\frac{1}{2} h_n$. Then, given
$t\in [0,1)$, we have, almost surely: \beq \limn \sup_{\wth_n\le
h\le h_n}\;\inf_{g\in \sqrt{2}\SSS}\Mmi
\frac{\DD_{n,h,t}}{(2h\log_2 n)^{1/2}}- g\Mmi =0.\label{vfau}\eeq
\end{theo}
The proof of Theorem \ref{T1} is written in \S \ref{trois}.\lb
\textbf{Remark}: Condition $\wth_n<h_n/2$ is just technical, as
this result is really interesting when $(\wth_n)\suite$ and
$(h_n)\suite$ are sequences that tend to 0 at different rates
 (typically $n^{-\alp_1}$ and $n^{-\alp_2},
 \;0<\alp_1<\alp_2<1$). Clearly, Theorem \ref{T1} seems unsatisfactory, as one would expect the limit set to be $\SSS$
 instead of $\sqrt{2}\SSS$. As it will be pointed out in the proof of Theorem \ref{T1} (see \S \ref{trois2}), it is possible to prove (\ref{fau}) when
 \beq \forall \beta>0,\;\limn \log(h_n/\wth_n)/(\log n)^{\beta}=0.\label{epqe}\eeq
However, (\ref{epqe}) is a very restrictive condition, imposing
$(h_n)\suite$ and $(\wth_n)\suite$ to have rates of convergence to
zero that are very close one to each other. In \S
 \ref{trois}, we shall try to point out the main
 difficulty that imposes us to weaken (\ref{fau}) to (\ref{vfau}). Showing that (\ref{fau}) is true or false without imposing (\ref{epqe}) remains an
 open problem.\lb
The second step of our investigation is to determine the validity
of (\ref{faux}). This assertion turns out to be false as soon as
$\wth_n/h_n\rar 0$, which is a consequence of our next result. We
first need to introduce some further notations. Given an integer
$k\geq 2$, we endow the space $(B[0,1])^k$ with the product
sup-norm, namely $\mmi g_1,\ldots,g_k\mmi_k:=\max\{\mmi g_1
\mmi,\ldots,\mmi g_k\mmi\}$, and we define \beq\SSS_k:=\aoo
(g_1,\ldots,g_k)\in (AC[0,1])^k, \sli_{j=1}^k\mmi g_j
\mmi_H^2\le1\aff\label{SSS2}. \eeq Now consider sequences
$0<h_{n,1}<\ldots<h_{n,k}<1$ satisfying, ultimately as $\nif$,
\begin{align}
h_{n,l}/h_{n,l+1} \downarrow 0,\;l=1,\ldots, k-1,\label{multi1}\\
h_{n,k}\downarrow 0,\; nh_{n,1}\uparrow \infty\label{multi2}.
\end{align}
Our second main result is the following functional limit law,
which is proved in \S \ref{quatre}.
\begin{theo}\label{T2}
Under assumptions (\ref{multi1}) and (\ref{multi2}) we have almost
surely \beq \poo \frac{\DD_{n,h_{n,1},t}}{(2h_{n,1}\log_2
n)^{1/2}},\ldots,\frac{\DD_{n,h_{n,k},t}}{(2h_{n,k}\log_2
n)^{1/2}}\pff\leadsto \SSS_k\label{Doubleloc}.\eeq Here $\leadsto$
refers to the Banach space $\po B([0,1])^k,\norm_k\pf$.\end{theo}
Note that $\SSS_k$ is the unit ball of the reproducing kernel
Hilbert space of $(W_1,\ldots,W_k)$, where $W_1,\ldots,W_k$ are
independent Wiener processes on $[0,1]$. Theorem \ref{T2}
describes an asymptotic independence phenomenon which has been
earlier investigated by Deheuvels (2000) and Deheuvels \textit{et
al.} (1999). The proof of Theorem \ref{T2} is provided in \S
\ref{quatre}. Now, to see that (\ref{faux}) is false, choose $g$
as the identity function so as $(g,g)\notin \SSS_2$, which entails
that $\inf\{\mmi g-g_1,g-g_2\mmi_2,\;(g_1,g_2)\in \SSS_2\}>\e_0$
for some $\e_0>0$. By Theorem \ref{T2} we have, almost surely,
\begin{align}
\nono &\lin \sup_{\wth_n\le h\le h_n}\Mmi
\frac{\DD_{n,h,t}}{(2h\log_2 n)^{1/2}}-g\Mmi\\
\nono\geq &\lin \Mmi \frac{\DD_{n,\wth_n,t}}{(2\wth_n \log_2
n)^{1/2}}-g,\;\frac{\DD_{n,h_n,t}}{(2h_n \log_2 n)^{1/2}}-g\Mmi_2\\
\nono\geq& \e_0,
\end{align}
which invalidates (\ref{faux}).\vk \textbf{A local Bahadur-Kiefer
representation}\lb A consequence of Theorem \ref{T1} is the
following local Bahadur-Kiefer representation, which is very
largely inspired from Einmahl and Mason (1988, Theorem 5). For
$0<h<1$ and $n\geq 1$ we set $a_n(h):= (h\log_2
n/n)^{1/2},\;b_n(h):=\log(nh),\;d_n(h):=2\log_2 n+b_n(h)$,
$r_n(h):=(a_n(h)d_n(h))^{1/2}$ and
$$R_n(h):= \Mmi
\DD_{n,h,0}+\DD'_{n,h,0}\Mmi.$$
\begin{theo}\label{T3}
Under the conditions of Theorem \ref{T1}, with $t=0$, we have,
almost surely \beq \lsn \sup_{\wth_n\le h\le
h_n}r_n(h)^{-1}R_n(h)\le 2^{1/2}.\label{bk}\eeq
\end{theo}
The proof of Theorem \ref{T3} is provided in \S \ref{cinq}.\lb
\textbf{Remark}: In view of Theorem 5 of Einmahl and Mason (1988),
Theorem \ref{T3} seems to be non optimal since a factor $2^{1/4}$
can be drop when $h_n=\wth_n$. This is a consequence of the fact
that we were only able to prove (\ref{vfau}) instead of
(\ref{fau}).
\section{Proof of Theorem \ref{T1}}\label{trois}
Our proof is divided into two subsections. In \S \ref{trois1}, we
establish a large deviation result which holds uniformly in
$\wth_n\le h\le  h_n$. Then we make use of that (uniform) large
deviation principle to prove Theorem \ref{T1} in \S \ref{trois2}.
\subsection{A uniform large deviation principle} \label{trois1}
\subsubsection{Definitions}
Large deviation results are commonly used when proving functional
laws of the iterated logarithm such as (\ref{Masonloc}). As a
uniformity in $\wth_n\le h\le h_n$  appears in Theorem \ref{T1},
we shall make use of a large deviation principle that holds
uniformly in $h$. This tool was first used by Mason (2004). From
now on, $(\e_{n,i})\tab$ will denote a triangular array of
strictly positive numbers satisfying $\max_{1\le i\le p_n}
\e_{n,i}\rar 0$ as $\nif$. We call a rate function in a metric
space $(E,d)$ any positive real function $J$ on $E$ such that, for
each $a\geq 0$, the set $\{g\in E,\;J(g)\le a\}$ is a compact set
of $(E,d)$.
\begin{defi}
Let $(E,d)$ be a metric space and let $\TT_0$ be a $\sig$-algebra
included in the Borel $\sig$-algebra of $(E,d)$. Let
$(X_{n,i})\tab$ be a triangular array of random variables that are
measurable for $(E,\TT_0)$. We say that $(X_{n,i})\tab$ satisfies
the uniform large deviation principle (ULDP) for $(\e_{n,i})\tab$,
a rate function $J$ and $\TT_0$ whenever
\begin{enumerate}
  \item For each closed set $F\in \TT_0$ we have \beq\lsn \max_{i\le p_n} \e_{n,i} \log \poo\PPP\poo X_{n,i}\in F\pff\pff\le -J(F),\label{GDU1}\eeq
  \item For each open set $O\in \TT_0$ we have \beq\lin \min_{i\le p_n} \e_{n,i} \log \poo\PPP\poo X_{n,i}\in O\pff\pff\geq
  -J(O).\label{GDU2}\eeq
\end{enumerate}
\end{defi}
\textbf{Remark:} In this definition, we introduce a sub
$\sig$-algebra $\TT_0$ because we will consider repeatedly $(E,d)$
as the metric space $(B([0,1],\norm)$. As the $\DD_{n,h,t}$ are
not Borel measurable in that space, we shall consider $\TT_0$ as
the $\sig$-algebra spawned by the open balls of $(B([0,1],\norm)$.
We will sometimes take $(E,d)$ as a finite dimensional vector
space, in which case $\TT_0$ will denote the Borel $\sig$-algebra.
Another way to avoid measurability problems is to consider inner
and outer probabilities (see, e.g.,Van der Vaart and Wellner
(1996), Chapter 1).\lb The next result is a consequence of the
work of Arcones (2003).
\begin{propo} \label{GDUstoch}
Let $(X_{n,i})\tab$ be a triangular array of random variables
taking values in $B([0,1])$ and measurable for $\TT_0$. Let
$(\e_{n,i})\tab$ be a triangular array of strictly positive real
numbers. Assume that the following conditions hold:
\begin{enumerate}
\item For each $p\geq 1$ and $(s_1,\ldots,s_p)\in (0,1)^p$
satisfying $s_i\neq s_j$ for each $i\neq j$, the triangular array
$\po X_{n,i}(s_1),\ldots,X_{n,i}(s_p)\pf\tab$ satisfies the ULDP
in $\RRR^p$ for $(\e_{n,i})\tab$ and a rate function
$I_{s_1,\ldots,s_p}$. \item For any $\tau>0$ we have
$$\lim_{\dd\downarrow 0}\lsn \max_{i\le p_n} \log\poo \PPP \poo \sup_{\mid
s-s' \mid <\dd} \mid
X_{n,i}(s')-X_{n,i}(s)\mid>\tau\pff\pff=-\infty.$$
\end{enumerate}
Then $(X_{n,i})\tab$ satisfies the ULDP in $(B([0,1],\norm)$ for
$(\e_{n,i})\tab$, $\TT_0$ and the following rate function:
$$I(g):=\sup_{p\geq 1,\;(s_1,\ldots,s_p)\in (0,1)^p}
I_{s_1,\ldots,s_p}\po g(s_1),\ldots,g(s_p)\pf,\;g\in B([0,1]).$$
\end{propo}
Now consider the following rate function on $B([0,1])$ that is
known to rule the large deviation properties of a Wiener process:
\beq J(g):=\left\{%
\begin{array}{ll}
    \mmi g\mmi_H^2, & \hbox{when $g\in AC[0,1]$;} \\
    \infty, & \hbox{when $g\notin AC[0,1]$.} \\
\end{array}%
\right.  \label{J}\eeq Notice that $\SSS=\{g\in
B([0,1]),\;g(0)=0,\;J(g)\le 1\}.$
 The main tool that will be used to
achieve our proof of Theorem \ref{T1} is the following ULDP.
\begin{propo}\label{GDU}
Let $(h_n)\suite$ and $(\wth_n)\suite$ be two sequences satisfying
conditions of Theorem \ref{T1} and let $(h_{n,i})\tab$ be a
triangular array satisfying $\wth_n \le h_{n,i} \le h_n$ for each
$n\geq 1,\;i\le p_n$. Then the triangular array
$$\poo(2h_{n,i}\log_2 n)^{-1/2} \DD_{n,h_{n,i},t}\pff\tab$$
satisfies the ULDP in $(B([0,1]),\norm)$ for $\TT_0$, the rate
function $J$ given in (\ref{J}) and the (constant in $i\le p_n$)
triangular array $(1/\log_2 n)\tab$.
\end{propo}
\textbf{Proof of Proposition \ref{GDU}}: We shall make use of
Proposition \ref{GDUstoch}, and we hence have to show that
conditions 1 and 2 of this proposition are satisfied. This
verification will be a consequence of two separate lemmas. The
next proposition, which shall be useful to prove our first lemma,
follows directly from the arguments of Ellis (1984). Here
$<\cdot,\cdot>$ denotes the usual scalar product on $\RRR^p$.
\begin{propo}\label{Ellis}
Let $(X_{n,i})\tab$ be a triangular array of random vectors taking
values in $\RRR^p$, and let $(\e_{n,i})\tab$ be a triangular array
of strictly positive real numbers. Assume that there exists a
positive real function $\ell$ (which may take infinite values) on
$\RRR^p$ such that the following conditions are satisfied.
\begin{enumerate}
\item $\ell$ is convex and lower semi continuous on $\RRR^p$.
\item The definition set $D(\ell):=\{\lab \in \RRR^p,\;\ell
(\lab)<\infty\}$ has an interior that contains the null vector.
\item $\ell$ is differentiable on the interior of $D(\ell)$ and,
for each sequence $(\lab_n)\suite$ converging to a boundary point
of $D(\ell)$ we have $\mmi \nabla \ell (\lab_n)\mmi_{\RRR^p}\rar
\infty.$ Here $\norm_{\RRR^p}$ denotes the usual Euclidian norm.
 \item For each $\lab\in D(\ell)$, we have
$$\limn \max_{i\le p_n} \Mid \e_{n,i}\log\poo \EEE\poo \exp\po
\e_{n,i}^{-1} <\lab,X_{n,i}>\pf\pff\pff- \ell(\lab)\Mid=0.$$ \item
For each $\lab\notin D(\ell)$, we have
$$\limn \min_{i\le p_n}  \e_{n,i}\log\poo \EEE\poo \exp\po
\e_{n,i}^{-1} <\lab,X_{n,i}>\pf\pff\pff=\infty.$$
\end{enumerate}
Then $(X_{n,i})\tab$ satisfies the ULDP in $\RRR^p$ for
$(\e_{n,i})\tab$ with the following rate function:
$$\wtJ(s):=\sup_{\lab \in \RRR^p} <\lab,s>- \ell(\lab),\;s\in
\RRR^p.$$
\end{propo}
We now state our first lemma.
\begin{lem}\label{lem1}
Let $p\geq 1$ and $(s_1,\ldots,s_p)\in [0,1]^p$ be arbitrary, with
$s_1<s_2<\ldots<s_p$. Under the assumptions of Proposition
\ref{GDU}, the triangular array of $\RRR^p$-valued random vectors
$$\poo(2 h_{n,i}\log_2 n)^{-1/2}
\po\DD_{n,h_{n,i},t}(s_1),\ldots,\DD_{n,h_{n,i},t}(s_p)\pf\pff\tab$$
satisfies the ULDP for $(\e_{n,i})\tab$ with the following rate
function (with $s_0:=0$).
$$J_{s_1,\ldots,s_p}(x_1,\ldots,x_p):=\sli_{i=0}^p (s_{i+1}-s_i){\poo
\frac{x_{i+1}-x_i}{s_{i+1}-s_{i}}\pff}^2,\;(x_1,\ldots,x_p)\in
\RRR^p.$$
\end{lem}
\begin{plem}\end{plem} We shall make use of Proposition \ref{Ellis}. Fix
$\lab=(\lab_1,\ldots,\lab_p) \in \RRR^p$ and and write the
$\DD_{n,h_{n,i},t}$ as sums of i.i.d. random variables, namely
\beq (2h_{n,i}\log_2 n)^{-1/2}\sli_{j=1}^p
\lab_j\DD_{n,h_{n,i},t}(s_j)=(2nh_{n,i}\log_2
n)^{-1/2}\sli_{k=1}^n Z_{n,h_{n,i},t}^k,\eeq where
$$Z_{n,h_{n,i},t}^k:=\sli_{j=1}^p\lab_j
\po1_{[t,t+h_{n,i}s_j]}(U_k)-h_{n,i}s_j\pf,\; k=1,\ldots,n.$$
These $n$ random variables are i.i.d with mean 0 and
variance-covariance matrix given by $h_{n,i}\lab' \Sig_{n,i}
\lab$, with $\Sig_{n,i}(l,l'):= \min(s_l,s_{l'})-h_{n,i} s_l
s_{l'}.$ Now define the matrix $\Sig(l,l'):= \min(s_l,s_{l'})$.
Clearly, as $h_{n,i}\le h_n\rar 0$ we have $\Sig_{n,i}\rar \Sig$
uniformly in $i$ as $\nif$. By standard computations we have, for
each $n\geq 1$ and $i\le p_n$:
\begin{align}
\nono&(\log_2 n)^{-1}\log\poooo \EEE\pooo \exp \poo \log_2 n (2
h_{n,i}\log_2 n)^{-1/2}
\sli_{j=1}^p \lab_j\DD_{n,h_{n,i},t}(s_j)\pff\pfff\pffff\\
=& \frac{n}{\log_2 n}\log \poo\EEE \exp \po r_{n,i}
Z_{n,h_{n,i},t}^1\pf\pff,\label{plo2}
\end{align}
where $r_{n,i}:= (\log_2 n/2nh_{n,i})^{1/2}$. Recall that
$\max_{i\le p_n} r_{n,i}\rar 0$ as $\nif$, since $\wth_n$
satisfies (\ref{CRS}), and notice that the $Z^k_{n,h_{n,i},t}$ are
centered and almost surely bounded by $p\max_{j=1,\ldots,p} \mid
\lab_j\mid$. This ensures that the following Taylor expansion is
valid, for each $n\geq 1,\;i\le p_n$ (here $\varepsilon$ denotes a
real function satisfying $\varepsilon(u)\rar 0$ as $u\rar 0$):
\beq \EEE\poo \exp \po r_{n,i} W_{n,h_{n,i},t}^1\pf\pff= 1+
\frac{r_{n,i}^2 h_{n,i}}{2}\lab' \Sig_{n,i} \lab(1+
\varepsilon(r_{n,i}))\label{plo1}.\eeq Combining (\ref{plo2}) and
(\ref{plo1}), we get \begin{align} \nono \limn &\max_{i\le p_n}
\gMid \frac{\log\pooo \EEE\pooo \exp \poo\frac{\log_2 n}{(2h_{n,i}
\log_2 n)^{1/2}} \sli_{j=1}^p
\lab_j\DD_{n,h_{n,i},t}(s_j)\pff\pfff\pfff}{\log_2 n}- \frac{1}{4}
\lab' \Sig \lab\gMid=0.\end{align} As the function $\ell(\lab):=
\lab' (\Sig/4) \lab$ obviously satisfies conditions of Proposition
\ref{Ellis}, the proof of Lemma \ref{lem1} is concluded by
noticing that
$$\sup_{t\in\RRR^p} <t,x>-\ell(t)= x' \Sig^{-1} x=\sli_{i=0}^p
(s_{i+1}-s_i){\poo \frac{x_{i+1}-x_i}{s_{i+1}-s_{i}}\pff}^2.\Box$$
Our next lemma shows that condition 2 of Proposition
\ref{GDUstoch} is fulfilled.
\begin{lem}\label{lem2}
Under the assumptions of Proposition \ref{GDU}, we have, for each
$\tau>0$\beq\nono \lim_{\dd\downarrow 0} \lsn \max_{i\le p_n}
\log\poo\PPP\poo \sup_{\mid s-s'\mid<\dd}\Mid
\frac{\DD_{n,h_{n,i},t}(s)-\DD_{n,h_{n,i},t}(s')}{(2h_{n,i}\log_2
n)^{1/2}}\Mid\geq \tau\pff\pff=-\infty.\eeq
\end{lem}
\begin{plem}\end{plem} Fix $\tau>0$ and introduce a parameter $\dd>0$ that will be
chosen small enough in the sequel. The proof of this lemma relies
on an exponential inequality for the oscillations of the local
empirical process, which is due to Einmahl and Mason (1988) (see
their Inequality 1). For positive numbers $a,b$ with $a+b\le 1$,
write \beq \omega_n(a,b):=\mathop{\sup_{0\le s\le b,}}_{0\le s'\le
a} \mid \alp_n(s+s')-\alp_n(s)\mid\label{omegan} .\eeq
\begin{fact}[Einmahl, Mason, 1988] \label{f1}Fix $0< \varepsilon \le 1/2$.
There exists $K(\varepsilon)<\infty$ such that, for any $n\geq 1$,
$\lab>0$, $a>0,\;b>0$ fulfilling $a+b\le 1$ and $0<a<1/4$, \beq
\PPP\poo \omega_n(a,b) \geq \lab\pff\le
K(\varepsilon)ba^{-1}\exp\poo
-\frac{(1-\varepsilon)\lab^2}{2a}\Psi\poo\frac{\lab}{\sqrt{n}
a}\pff\pff.\label{acclocal}\eeq \end{fact}Here we write
$\Psi(u):=2u^{-2} ((1+u)\log(1+u)-u)$.\lb Applying
(\ref{acclocal}) to $b=h_{n,i},\;a=\dd h_{n,i},\;\varepsilon=1/2$
and $\lab=\tau(2h_{n,i}\log_2 n)^{1/2}$ we get, for all large $n$
and $i\le p_n$ (so that $h_{n,i}\le h_n\le 1/4$)
\begin{align}
\nono \PPP\poo \sup_{\mid s-s'\mid<\dd}\Mid
\frac{\DD_{n,h_{n,i},t}(s)-\DD_{n,h_{n,i},t}(s')}{(2h_{n,i}\log_2
n)^{1/2}}\Mid\geq \tau\pff\le&\frac{K(\frac{1}{2})}{\dd}\exp\pooo
-\frac{\tau^2 \log_2 n}{2\dd}\Psi\poo
\frac{\tau\sqrt{2\log_2n}}{\dd\sqrt{nh_{n,i}}}\pff\pfff\\
\le&\frac{K(\frac{1}{2})}{\dd}\exp\poo -\frac{\tau^2 \log_2
n}{4\dd}\pff\label{poinuh}.
\end{align}
The last inequality holds for all large $n$ and $i\le p_n$ since
$\Psi(u)\rar 1$ as $u\rar 0$, and since\beq \limn \max_{i\le p_n}
\frac{\log_2 n}{nh_{n,i}}=0.\label{tt}\eeq Now taking the
logarithm in (\ref{poinuh}) concludes the proof of Lemma
\ref{lem2}, then lemmas \ref{lem1} and \ref{lem2} in combination
with Proposition \ref{Ellis} conclude the proof of Proposition
\ref{GDU}. $\Box$
\subsection{Proof of Theorem
\ref{T1}}\label{trois2} We shall invoke usual blocking arguments
along the following subsequence: \beq n_k:=\coo\exp \poo
k\exp\po-(\log k)^{1/2}\pf\pff\cff,\;k\geq 5.\label{nk}\eeq
Clearly, $n_k$ satisfies, as $\kif$, \beq \frac{n_k}{n_{k+1}}\rar
1,\;\log_2(n_k)=\log k(1+o(1)).\label{propnk}\eeq Now define the
blocks $ N_k:=\{n_{k-1},\ldots,n_k-1\}$ for $k\geq 6.$  Fix $\e>0$
and consider a parameter $\rho>1$ that will be chosen small enough
in the sequel. For any $k\geq 5$, consider the following
discretisation of $[\wth_{n_k},h_{n_{k-1}}]$
\begin{align}
 h_{n_k,R_k}:=h_{n_{k-1}},\;\;
\hkl:=&\rho^l \wth_{n_k},\;l=0,\ldots, R_k-1,\label{discret1}
\end{align}
where $R_k:= [(\log (h_{n_{k-1}}/\wth_{n_k}))/\log(\rho)]+1$, and
$[u]$ denotes the only integer $q$ fulfilling $q\le u<q+1$.
Clearly, as $\kif$, we have\beq R_k=O(\log n_k)\label{grando}.\eeq
Our aim is to show that the following probabilities are summable
in $k$ so as the Borel-Cantelli lemma would complete the proof of
Theorem \ref{T1}. \beq \PPP_k:= \PPP\pooo \max_{n\in
N_k}\;\sup_{\wth_n\le h \le h_n} \inf_{g\in \sqrt{2}\SSS}\Mmi
\frac{\DD_{n,h,t}}{(2h\log_2
n)^{1/2}}-g\Mmi\geq3\e\pfff.\label{PPPk}\eeq  Clearly we have
\begin{align}
\nono\PPP_k\le&\;\PPP\pooo \max_{0\le l \le R_k}\inf_{g\in
\sqrt{2}\SSS}\Mmi
\frac{\DD_{n_k,\hkl,t}}{(2\hkl\log_2 \nk)^{1/2}}-g\Mmi\geq\e\pfff\\
\nono&\;+ \PPP\pooo \max_{n\in N_k}\;\max_{0\le l\le
R_k-1}\;\sup_{\hkl\le h\le \rho \hkl}\Mmi
\frac{\DD_{n,h,t}}{(2h\log_2
n)^{1/2}}-\frac{\DD_{n_k,\hkl,t}}{(2\hkl \log_2
n_k)^{1/2}}\Mmi>2\e\pfff\\
\nono =:& \;\PPP_{1,k}+\PPP_{2,k}.
\end{align}
To show that $\PPP_{1,k}$ is summable, we shall make use of
Proposition \ref{GDU}. Consider the following subset of
$B([0,1])$:
$$F:=\aoo f\in B([0,1]),\;\inf_{g\in \sqrt{2}\SSS} \mmi f-g\mmi \geq
\e\aff.$$ Since the rate function $J$ given in (\ref{J}) is lower
semi continuous on $(B([0,1],\norm)$, there exists $\alp_1>0$
satisfying $J(F)=2+2\alp_1$. Hence, for all large $k$ we have \beq
\PPP_{1,k}\le(R_k+1) \exp\po-(2+\alp_1)\log_2
n_k\pf.\label{ret1}\eeq Recalling (\ref{propnk}) and
(\ref{grando}), we conclude that $\PPP_{1,k}$ is summable in $k$.
It remains to show the summability of  $(\PPP_{2,k})_{k\geq 1}$.
First notice that
\begin{align}
\nono\PPP_{2,k}\le&\; \PPP\pooo \max_{l\le R_k-1}\max_{n\in
N_k}\;\sup_{\hkl\le h\le \rho\hkl}\Mmi
\frac{\sqrt{n}\DD_{n,h,t}-\sqrt{n}\DD_{n,h_{n_k,l},t}}{(2n_k\hkl\log_2
n_k)^{1/2}}\Mmi>\e\pfff\\
\nono&\;+\PPP\pooo\max_{l\le R_k-1}\max_{n\in N_k}\;\sup_{\hkl\le
h\le \rho\hkl}\BB(n,h)\Mmi \frac{\sqrt{n}\DD_{n,h,t}}{(2n_k\rho
\hkl\log_2
n_k)^{1/2}}\Mmi>\e\pfff\\
=:&\;\PPP_{3,k}+\PPP_{4,k}\label{baf},
\end{align}
where \beq \BB(n,h):=\Mid\sqrt{\frac{n_k \rho\hkl \log_2
n_k}{nh\log_2 n}}-1\Mid,\;n\in N_k,\;l\le R_k-1,\;\hkl\le h\le
\rho\hkl. \label{Bnh}\eeq We shall require a maximal inequality
due to Montgomery-Smith (1993) (see also Latala (1993)).
\begin{fact}[Montgomery-Smith, Latala, 1993] There exists a constant $c>0$ such that, given a Banach space
$(E,\norm)$ and a finite sequence $(X_i)_{1\le i\le n}$ of i.i.d.
random variables taking values in $(E,d)$ we have, for each
$\lab>0$: \beq\PPP\pooo \max_{1\le i\le n} \Mmi \sli_{j=1}^i X_j
\Mmi \geq \lab\pfff\le c\PPP\pooo \Mmi\sliin X_i\Mmi\geq
\frac{\lab}{c}\pfff.\label{Montgom}\eeq
\end{fact}
Applying inequality (\ref{Montgom}), we get
\begin{align}
\nono\PPP_{3,k}\le& \sli_{l=0}^{R_k-1}\PPP\pooo\max_{n\in
N_k}\;\sup_{\hkl\le h\le \rho\hkl}\Mmi
\frac{\sqrt{n}\DD_{n,h,t}-\sqrt{n}\DD_{n,h_{n_k,l},t}}{(2n_k\hkl\log_2
n_k)^{1/2}}\Mmi>\e\pfff\\
\le & c\sli_{l=0}^{R_k-1}\PPP\pooo\sup_{\hkl\le h\le
\rho\hkl}\Mmi\frac{\sqrt{n_k}\DD_{n_k,h,t}-\sqrt{n_k}\DD_{n_k,h_{n_k,l},t}}{(2n_k\hkl\log_2
n_k)^{1/2}}\Mmi>\e/c\pfff\label{thi}.
\end{align}
As $\hkl\le h_{\nkm}\rar 0,$ each term of (\ref{thi}) can be
bounded by inequality (\ref{acclocal}), provided that
$h_{\nkm}<1/4$. In inequality (\ref{acclocal}), we repeatedly
choose $b= \hkl,\;a=\hkl(\rho-1),\;\varepsilon=1/2,\;\lab=(2\hkl
\log_2 n_k)^{1/2}\e/c$. Hence, for all large $k$ we have
\begin{align}\nono\PPP_{3,k}\le&
c\sli_{l=0}^{R_k-1}\frac{K(\frac{1}{2})}{\rho-1}\exp \pooo
-\frac{\e^2\log_2 n_k}{2c^2(\rho-1)^2}\Psi\poo \frac{\e
\sqrt{\log_2 n_k}}{c(\rho-1)\sqrt{n_k \hkl}}\pff\pfff\\
\le &c\sli_{l=0}^{R_k-1}\frac{K(\frac{1}{2})}{\rho-1}\exp \pooo
-\frac{\e^2\log_2 n_k}{4c^2(\rho-1)^2}\pfff\label{gfpr}\\
\le&\frac{cK(\frac{1}{2})}{\rho-1}R_k
{k}^{-{\e/2c(\rho-1)}^2}\label{frfr}.
\end{align}
Inequality (\ref{gfpr}) is true for all large $k$ since
$\Psi(u)\rar 1$ as $u\rar 0$, and since \beq \limk \max_{l\le
R_k-1} \frac{\log_2 n_k}{n_k \hkl}=0\label{unifk}.\eeq Inequality
(\ref{frfr}) takes in account the fact that $\log_2 n_k= \log
k(1+o(1))$ as $\kif$. Hence for any choice of
$1<\rho<1+\sqrt{\e/2c}$ the general term (\ref{frfr}) is summable
in $k$ and so are the $\PPP_{3,k}$ (recall (\ref{grando})).
Showing that $\sum\PPP_{4,k}<\infty$ will be done in a similar
way. First notice that, as $n_k/\nkm \rar 1$ and $1\le \rho
\hkl/h\le \rho$ we have \beq \limk \max_{0\le l \le
R_k-1}\;\max_{n \in N_k} \BB(h,n)=\rho^{1/2}-1\le
2(\rho-1).\label{limrho}\eeq Hence, for all large $k$ we have
\begin{align}
\nono \PPP_{4,k}\le &\;\PPP\pooo \max_{0\le l\le
R_k-1}\;\max_{n\in N_k}
\;\Mmi\frac{\sqrt{n}\DD_{n,\rho\hkl,t}}{(2n_k
\rho\hkl \log_2 n_k)^{1/2}}\Mmi>\frac{\e}{2(\rho-1)}\pfff\\
\nono \le&\; c \sli_{l=0}^{R_k-1} \PPP\pooo
\Mmi\frac{\DD_{n_k,\rho\hkl,t}}{(2\rho\hkl \log_2
n_k)^{1/2}}\Mmi>\frac{\e}{2c(\rho-1)}\pfff\\
\le &\; 2c\sli_{l=0}^{R_k-1}\exp\pooo -
\frac{\e^2(1-\rho\hkl)\log_2
n_k}{8c^2(\rho-1)^2}\Psi\poo\frac{\e(1-\rho\hkl)\sqrt{2 \log_2
n_k}}{2c\sqrt{n_k \rho \hkl}}\pff\pfff\label{v1}\\
\le&\; 2c R_k \exp\pooo -\frac{\e^2(1-\rho\hkl)\log_2
n_k}{16c^2(\rho-1)^2}\pfff.\label{v2}
\end{align}
Here, (\ref{v1}) is a consequence of Inequality 2 in Shorack and
Wellner (1986, p. 444), with $p=\rho
\hkl,\;\lab=\e(1-\rho\hkl)(2\rho \hkl \log_2
n_k)^{1/2}/4c(\rho-1)$. Recalling (\ref{unifk}), we see that
(\ref{v2}) holds for all large $k$, as $\Psi(u)\rar 1$ when $u\rar
0$. Now choosing $\rho>1$ small enough leads to he summability of
$(\PPP_{4,k})_{k\geq 1}$, which concludes the proof of Theorem
\ref{T1}.$\Box$ \lb \textbf{Remark}: If we had replaced the limit
set $\sqrt{2}\SSS$ by $\SSS$ in Theorem \ref{T1}, then
(\ref{ret1}) would become $$\PPP_{1,k}\le
(R_k+1)\exp\po-(1+\alp_1)\log_2 n_k\pf.$$ Hence, we would be able
to conclude that $\PPP_{1,k}$ is summable if the cardinality
$R_k+1$ of the grids were smaller than $(\log n_k)^\beta$ for any
$\beta>0$. When constructing the $\hkl$ as in (\ref{discret1}),
the just mentioned condition is violated as soon as $\wth_n$ and
$h_n$ have "really" different rates of convergence to zero
(typically when $\wth_n=h^{-\beta_1}<n^{-\beta_2}$ with
$0<\beta_2<\beta_1<1$). It seems however impossible to reduce the
cardinality $R_k+1$ of our grids, since the oscillations between
two consecutive $\hkl$ become hardly controllable and hence the
corresponding probabilities $\PPP_{2,k}$ might not be summable.
One could expect some improvements of this proof, since the RHS of
(\ref{ret1}) is crudely obtained, but this turns out to be non
trivial, as Proposition \ref{GDU} would have to be improved to
more accurate large deviation rates for the
$\DD_{n_k,\hkl,t},\;0\le l\le R_k$. Another possibility would be
to "poissonize"  the $\DD_{n,h,t}$ and then make use of strong
approximation of a centred Poisson process by a Wiener process $W$
(see $\mathrm{Koml\grave{o}s}$ \textit{et al.}, 1977), which would
reduce the problem to studying the summability of \beq
\PPP_{1,k}^W:=\PPP\poo \exists \rho \in
(\frac{\wth_{n_k}}{h_{\nkm}},1),\; \rho^{-1/2}W(\rho\cdot) \notin
(2\log_2 n_k)^{1/2}(\SSS+\e\BB_0)\label{grpae}\pff,\eeq and then
try to make use of the isoperimetric properties of a Gaussian
measures (here $\BB_0$ denotes the unit ball of $B([0,1])$). This
however fails to work by making brute use of the isoperimetric
inequality, as long as $\wth_{\nk}/h_{\nkm}$ is not negligible
with respect to $\log_2 n_k$ as $\kif$. We hope however, that
(\ref{grpae}) may be better controlled and we thus leave an open
question to specialists in Gaussian measures.\lb
\section{Proof of
Theorem \ref{T2}}\label{quatre} To avoid lengthy notations, we
shall prove Theorem \ref{T2} only with $k=2$ with no loss of
generality. The key of our proof of Theorem \ref{T2} is the
following lemma.
\begin{lem}\label{lem3}
Under the assumptions of Theorem \ref{T2}, for any $p\geq 1$, $0<
s_1^{(1)}<\ldots<s_p^{(1)}<1$ and $0<
s_1^{(2)}<\ldots<s_p^{(2)}<1$, the sequence of $\RRR^{2p}$-valued
random vectors
$$X_n:=\pooo
\frac{\DD_{n,h_{n,1},t}(s_1^{(1)})}{(2 h_{n,1} \log_2
n)^{1/2}},\ldots, \frac{\DD_{n,h_{n,1},t}(s_p^{(1)})}{(2 h_{n,1}
\log_2 n)^{1/2}},\frac{\DD_{n,h_{n,2},t}(s_1^{(2)})}{(2 h_{n,2}
\log_2 n)^{1/2}},\ldots,\frac{\DD_{n,h_{n,2},t}(s_p^{(2)})}{(2
h_{n,2} \log_2 n)^{1/2}}\pfff$$ satisfies the large deviation
principle for the sequence $(\log_2 n)^{-1}$ and the following
rate function (writing $s^{(1)}_0=s^{(2)}_0=0$). \begin{align}
\nono
\ovJ_{s_1^{(1)},\ldots,s_p^{(1)},s_1^{(2)},\ldots,s_p^{(2)}}(x):=&\sli_{i=1}^p
(s^{(1)}_{i+1}-s_i^{(1)}){\poo\frac{x^{(1)}_{i+1}-x_i^{(1)}}{s_{i+1}^{(1)}-s_i^{(1)}}\pff}^{2}+(s^{(2)}_{i+1}-s_i^{(2)}){\poo\frac{x^{(2)}_{i+1}-x_i^{(2)}}{s_{i+1}^{(2)}-s_i^{(2)}}\pff}^{2},\\
&x=x_1^{(1)},\ldots,x_p^{(1)},x_1^{(2)},\ldots,x_p^{(2)}\in
(0,1)^{2p}.
\end{align}
\end{lem}
\begin{plem}\end{plem} The proof follows the same lines as the proof of Lemma
\ref{lem1}. Choose
$\lab:=(\lab_1^{(1)},\ldots,\lab_p^{(1)},\lab_1^{(2)},\ldots,\lab_p^{(2)})\in
\RRR^{2p}$ arbitrarily and set (recall that $U_1$ is uniform on
$[0,1]$).
\begin{align}
\nono X_{n,1}:=\sli_{j=1}^p \lab_j^{(1)}\po1_{[t,t+h_{n,1}
s_j^{(1)}]}(U_1)-h_{n,1}
s_j^{(1)}\pf,\\
\nono X_{n,2}:=\sli_{j=1}^p \lab_j^{(2)}\po1_{[t,t+h_{n,2}
s_j^{(2)}]}(U_1)-h_{n,2} s_j^{(2)}\pf.
\end{align}
By independence we have \begin{align} \nono&(\log_2 n)^{-1}
\log\poo \EEE
\poo \exp \po \log_2 n\;<\lab,X_n>\pf\pff\pff\\
\nono=&\frac{n}{\log_2 n}\log \poo \EEE \poo\exp \po
r_{n,1}X_{n,1}+r_{n,2}X_{n,2}\pf\pff\pff,
\end{align}
with $r_{n,1}:=\sqrt{\log_2 n/2n h_{n,1}}$ and $r_{n,2}:=
\sqrt{\log_2 n/2nh_{n,2}}$. As $X_{n,1}$ (resp $X_{n,2}$) is
centered and almost surely bounded by $2p \max_{j=1,\ldots,
2p}\mid \lab_j\mid$, the following Taylor expansion is valid by
the dominated convergence theorem (here $\displaystyle{\lim_{ \mid
a\mid, \mid b\mid\rar 0}\varepsilon(a,b)=0}$):
\begin{align} \nono& \log \poo \EEE\poo \exp \po
r_{n,1}X_{n,1}+r_{n,2}X_{n,2}\pf\pff\pff\\
\nono =& \frac{1}{2}\poo r_{n,1}^{2} \Var(X_{n,1})+r_{n,2}^{2}
\Var(X_{n,2})+2r_{n,1}r_{n,2}\Cov(X_{n,1},X_{n,2})\pff(1+\varepsilon(r_{n,1},r_{n,2})).
\end{align}
Now, writing $\lab_1:=(\lab_1^{(1)},\ldots,\lab_p^{(1)})$ and
$\lab_2:=(\lab_1^{(2)},\ldots,\lab_p^{(2)})$ we can write
$\Var(X_{n,1})=\lab_1'\Sig^{(1)}_n \lab_1$ and
$\Var(X_{n,2})=\lab_2'\Sig_n^{(2)} \lab_2$, where
\begin{align}
\nono \Sig^{(1)}_n(i,j):=&
h_{n,1}\min(s_i^{(1)},s_j^{(1)})-h_{n,1}^{2}s_i^{(1)}s_j^{(1)},\;\mathrm{and}\\
\nono\Sig_n^{(2)}(i,j):=&h_{n,2}\min(s_i^{(2)},s_j^{(2)})-h_{n,2}^{2}
s_i^{(2)}s_j^{(2)}. \end{align} Hence, setting
$$\Sig^{(1)}(i,j):=\min(s_i^{(1)},s_j^{(1)})\;\;\mathrm{and}\;\;
\Sig^{(2)}(i,j):=\min(s_i^{(2)},s_j^{(2)}),$$ we obtain \beq \po
r_{n,1}^{2} \Var(X_{n,1})+r_{n,2}^{2} \Var(X_{n,2})\pf=
\frac{\log_ 2
n}{2n}\po\lab_1'\Sig^{(1)}\lab_1+\lab_2'\Sig^{(2)}\lab_2\pf(1+o(1))\label{fav1}.\eeq
In a similar way, we can write $\Cov(X_{n,1},X_{n,2})=
\lab_1'\Sig_n\lab_2$, where $\Sig_n(i,j):=\min(h_{n,1}
s_i^{(1)},h_{n,2} s_j^{(2)})- h_{n,1} h_{n,2} s_i^{(1)}s_j^{(2)}.$
Now recalling that $h_{n,1}/h_{n,2}\rar 0$ we have
$\Sig_n(i,j)=h_{n,1} s_i^{(1)}(1-s_j^{(2)} h_{n,2})$ for all large
$n$, whence \beq  \Mid
r_{n,1}r_{n,2}\Cov(X_{n,1},X_{n,2})\Mid=\frac{\log_2
n}{n}\sqrt{\frac{h_{n,1}}{h_{n,2}}}(1+o(1))=o\poo\frac{\log_2
n}{n}\pff.\label{fav2}\eeq Combining (\ref{fav1}) and (\ref{fav2})
we get \beq \nono\limn(\log_2 n)^{-1} \log\poo \EEE \poo \exp \po
\log_2 n\;<\lab,X_n>\pf\pff\pff= \frac{1}{4} (\lab_1'\Sig^{(1)}
\lab_1+ \lab_2'\Sig^{(2)} \lab_2).\eeq Then applying Proposition
\ref{Ellis} leads to the claimed result.$\Box$\lb We shall now
show that Lemma \ref{lem3} is sufficient to infer a large
deviation principle for the couples of processes $(2h_{n,1} \log_2
n)^{-1/2} \DD_{n,h_{n,1},t}$ and $(2h_{n,2} \log_2 n)^{-1/2}
\DD_{n,h_{n,2},t}$. Consider the following processes on $[0,2]$
that are obtained by concatenation of $(2h_{n,1} \log_2 n)^{-1/2}
\DD_{n,h_{n,1},t}$ with $(2h_{n,2} \log_2 n)^{-1/2}
\DD_{n,h_{n,2},t}$:
$$
\wt{\DD_n}(s):=\left\{%
\begin{array}{ll}
    \frac{\DD_{n,h_{n,1},t}(s)}{(2h_{n,1}\log_2
n)^{1/2}},&\mathrm{when\;0\le s\le 1}; \\
    \frac{\DD_{n,h_{n,2},t}(s-1)}{(2h_{n,2}\log_2
n)^{1/2}}, & \mathrm{when\;1< s\le 2}. \\
\end{array}%
\right.$$Combining Lemma \ref{lem3} with Lemma \ref{lem2} we
conclude that conditions of Proposition \ref{GDUstoch} are
fulfilled, and thus $\wt{\DD}_n$ satisfies the large deviation
principle for $\e_n:= (\log_2n)^{-1}$ and for the following rate
function:
\begin{align}
\nono& \ovJ(g)\\
\nono:= &\sup\aoo\sli_{j=0}^p (s_{j+1}^{(1)}-s_j^{(1)}){\poo
\frac{g(s_{j+1}^{(1)})-g(s_j^{(1)})}{s_{j+1}^{(1)}-s_j^{(1)}}\pff}^{2}+(s_{j+1}^{(2)}-s_j^{(2)}){\poo
\frac{g(1+s_{j+1}^{(2)})-g(1+s_j^{(2)})}{s_{j+1}^{(2)}-s_j^{(2)}}\pff}^{2}, \\
\nono&p\geq
1,\;0<s_1^{(1)}<\ldots<s_p^{(1)}<1<1+s_1^{(2)}<\ldots<1+s_p^{(2)}<2\aff\\
\nono=& \mmi g^{(1)} \mmi_H^{2}+\mmi g^{(2)} \mmi_H^{(2)},
\end{align}
where $g^{(1)}(s):= g(s),\;g^{(2)}(s):=g(1+s)$, $s\in [0,1]$. The
remainder of the proof of Theorem \ref{T2} is a routine use of
usual techniques in local empirical processes theory (refer, e.g.,
to Deheuvels and Mason (1990)). We omit details for sake of
briefness. $\Box$
\section{Proof of Theorem \ref{T3}}\label{cinq}
We shall proceed in three steps. Recall that $a_n(h):= (h\log_2
n/n)^{1/2},\;b_n(h):=\log(nh),\;d_n(h):=2\log_2 n+b_n(h)$,
$r_n(h):=(a_n(h)d_n(h))^{1/2}$ and $R_n(h):= \Mmi
\DD_{n,h,0}+\DD'_{n,h,0}\Mmi.$
\begin{lem}\label{lem21}
Under the assumptions of Theorem \ref{T1}, we have almost surely
\beq \lsn\sup_{\wth_n\le h\le h_n}\frac{\mmi
F_n^{\lar}(h\cdot)\mmi}{h}=1\label{sup1}.\eeq
\end{lem}
\begin{plem}\end{plem} First notice that, almost surely, for each
$\rho>1,\;h>0,\;n\geq 1$, \begin{align} \nono F_n^{\lar}(h)\le
\rho h&\Leftarrow \frac{\DD_{n,\rho h,0}}{(2 h \log_2
n)^{1/2}}+(\rho-1){\poo\frac{nh}{2 \log_2 n}\pff}^{1/2}\geq
0.\end{align} Now, for fixed $\rho>1$ we have $(\rho-1)\inf\{
nh/\log_2n,\;\wth_n\le h\le h_n\}\rar \infty$. Moreover, by a
straightforward use of Theorem \ref{T1} and (\ref{Masonloc}), \beq
\lin \inf_{\wth_n\le h \le h_n} \frac{\DD_{n,\rho h,0}}{(2 h
\log_2 n)^{1/2}}\geq-(2\rho)^{1/2}\text{ almost surely.}\eeq This
shows that $(\ref{sup1})$ holds with $\le$ instead of $=$, while
the converse inequality trivially holds by Kiefer (1972), Theorem
6.$\Box$\lb
\begin{lem}\label{lem22}
Under the assumptions of Theorem \ref{T1} we have almost surely $$
\lsn \sup_{\wth_n \le h \le h_n} \frac{ \mmi
\DD'_{n,h,0}\mmi}{(2h\log_2 n)^{1/2}}\le 2^{1/2}.$$\end{lem}
\begin{plem}\end{plem}From Inequality (2.23) in Einmahl and Mason (1988) we have,
for each $n\geq 1$ and $h>0$,
$$ \frac{\mmi\DD'_{n,h,0}\mmi}{(2 h \log_2 n)^{1/2}}\le
\frac{\mmi\DD_{n,F_n^{\lar}(h),0}\mmi}{(2 h \log_2
n)^{1/2}}+\frac{1}{(2n h \log_2 n)^{1/2}}.$$The second term can be
drop since $n\wth_n \rar \infty$. Fix $\rho>0$. By Lemma
\ref{lem21} we have almost surely, for all large $n$ and for all
$\wth_n\le h\le h_n$, $$ \frac{\mmi\DD_{n,F_n^{\lar}(h),0}\mmi}{(2
h \log_2 n)^{1/2}}\le \rho^{1/2} \frac{\mmi\DD_{n,\rho
h,0}\mmi}{(2 \rho h \log_2 n)^{1/2}}, $$ from where we readily
obtain, by Theorem \ref{T1},
$$ \lsn\sup_{\wth_n \le h \le h_n} \frac{n^{-1/2} \mmi
\DD'_{n,h,0}\mmi}{(2h\log_2 n)^{1/2}}\le
(2\rho)^{1/2}\;\mathrm{almost\;surely.}$$ As $\rho>1$ was
arbitrary, Lemma \ref{lem22} is proved. $\Box$\lb The expression
$\omega_n$ appearing in the next lemma has been defined in
(\ref{omegan}).
\begin{lem}\label{lem23}
Under the assumptions of Theorem \ref{T1}, and given $\eta>0$, we
have almost surely \beq \lsn \sup_{\wth_n\le h \le h_n}
\frac{\omega_n(\eta a_n(h),h)}{r_n(h)}\le
\eta^{1/2}.\label{limsup}\eeq
\end{lem}
\begin{plem}\end{plem}This proof is largely inspired from the proof of Lemma 6 in
Einmahl and Mason (1988). Fix $\e>0$ and consider the sequence
$(n_k)$ the sets $N_k$ and the grids $\hkl,\:0\le l\le R_k$ as in
\S \ref{trois2}. Also define, for each $k\geq 5$ and $l\le R_k$,
\begin{align}
\nono a_{k,l}:=&\eta(\rho\hkl\log_2\nk/\nkm)^{1/2}\;\;\text{and}\\
\nono r_{k,l}:=&(a_{k,l}(2\log_2\nk+\log (\nk \hkl)))^{1/2}.
\end{align} As $a_{k,l}\geq a_n(h)$ for
each $n\in N_k$ and $h\in [\hkl,\rho\hkl]$, we have
\begin{align}
\nono &\PPP\pooo \bculi_{n\in N_k}\bculi_{\wth_n\le h\le h_n}
\frac{\omega_n(\eta a_n(h),h)}{r_n(h)}\geq \eta^{1/2}(1+3\e)\pfff\\
\nono \le &\;\PPP\pooo \bculi_{l=0}^{R_k-1}\bculi_{n\in
N_k}\bculi_{\hkl\le h\le \rho\hkl}
\frac{\omega_n(a_{k,l},\rho\hkl)}{r_n(h)}>\eta^{1/2}(1+3\e)\pfff\\
\le & \PPP\pooo  \bculi_{l=0}^{R_k-1}\bculi_{n\in N_k}
\frac{\omega_n(a_{k,l},\rho\hkl)}{r_{k,l}}>\eta^{1/2}(1+2\e)\pfff\label{razop}\\
=:&\ovPPP_k,
\end{align}
where (\ref{razop}) holds for any choice of $\rho>1$ small enough,
ultimately as $\kif$, which is a consequence of the easily checked
fact that \beq \lim_{\rho\rar 1}\limk \max_{n\in N_k}\max_{l\le
R_k-1}\;\sup_{h\in[\hkl,\rho\hkl]} \Mid
\frac{r_n(h)}{r_{k,l}}-1\Mid=0.\eeq By Bonferroni's inequality we
can write
\begin{align} \nono\ovPPP_{k}\le&
\sli_{l=0}^{R_k-1}\PPP\pooo\bculi_{n\in N_k}
\frac{\omega_n(a_{k,l},\rho\hkl)}{r_{k,l}}>\eta^{1/2}(1+2\e)\pfff\\
=:&\sli_{l=0}^{R_k-1}\ovPPP_{k,l}\label{defo}.\end{align} Some
straightforward verifications show that the blocking arguments of
Inequality 2 in Einmahl and Mason (1988) can be used
simultaneously to each $\ovPPP_{k,l}$, for all large $k$ and
hence, by Fact \ref{f1},
\begin{align}
\nono\ovPPP_{k,l}\le&2 \PPP\poo \omega_{\nk}(a_{k,l},\rho\hkl)\geq
\eta^{1/2}r_{k,l}(1+\e)\pff\\
\nono\le & 2K\po\frac{\e}{2}\pf\frac{\rho\hkl}{a_{k,l}}\exp
\poo-\frac{(1-\frac{\e}{2})(1+\e)^2}{2a_{k,l}}\eta
r^2_{k,l}\Psi\po\Delta_{k,l}\pf\pff,
\end{align}
where
$\Delta_{k,l}:=(1+\e)\eta^{1/2}r_{k,l}\nk^{-1/2}a_{k,l}^{-1}$
converge to 0 uniformly in $l\le R_k-1$ when $\kif$. Since $\Psi$
(given in Fact \ref{f1}) satisfies $\Psi(u)\rar 1$ as $u\rar 0$ we
obtain, for all large $k$ and for each $l\le R_k-1$,
\begin{align}\nono \ovPPP_{k,l}\le& 2K\po\frac{\e}{2}\pf\sqrt{\frac{\nkm\rho\hkl}{\eta^2\log_2 \nk}}\exp\pooo -\frac{(1-\frac{\e}{2})^2(1+\e)^2}{2}\po 2\log_2 \nk+\log(\nk\rho\hkl)\pf\pfff\\
\nono\le&
2K\po\frac{\e}{2}\pf{\po\frac{\eta^2}{\rho}\pf}^{\e/8}(\nkm
\hkl)^{-\e/8}(\log_2 n_k)^{-1/2}(\log\nkm)^{-1-\e/4},\end{align}
for all large $k$ and for each $0\le l\le R_k-1$, which entails by
(\ref{defo})
\begin{align} \nono \ovPPP_{1,k}\le
&2K\po\frac{\e}{2}\pf{\po\frac{\eta^2}{\rho}\pf}^{\e/8}(\log_2\nk)^{-1/2}
(\log\nkm)^{-1-\e/4}\nkm^{-\e/8}h_{\nk}^{-\e/8}\sli_{l=0}^{R_k-1}\rho^{-l\e/8}\\
\nono\le
&2K\po\frac{\e}{2}\pf{\po\frac{\eta^2}{\rho}\pf}^{\e/8}\frac{1}{1-\rho^{-\e/8}}(\log_2\nk)^{-1/2}(\log\nkm)^{-1-\e/8}(\nkm
h_{\nk})^{-\e/8},\end{align} from where $\ovPPP_{k}$ is summable
in $k$. $\Box$ \lb The proof of Theorem \ref{T3} is concluded as
follows. First, it is well known that, almost surely, \beq
\mmi\alp_n+\beta_n+(\alp_n(F_n^{\lar})-\alp_n)\mmi=n^{-1/2},\label{ersq1}\eeq
whence, almost surely, for all $n\geq 1$ and $h>0$, \beq
R_n(h)\le\sup_{0<s<h}
\mmi\alp_{n}(s+n^{-1/2}\beta_n(s))-\alp_n(s)\mmi+n^{-1/2},\eeq
from where \beq \nono r_n(h)^{-1}R_n(h)\le r_n(h)^{-1}
\omega_n(n^{-1/2}\mmi\DD_{n,h,0}\mmi,h)+(nh \log_2
n)^{-1/4}(2\log_2 n+ \log(nh))^{-1/2},\label{tto}\eeq which
concludes the proof by combining lemmas \ref{lem22} and
\ref{lem23} (with the choice of $\eta=2$), as the second term of
the RHS of \ref{tto} converges to 0 uniformly in $\wth_n\le h\le
h_n$ as $\nif$.$\Box$

\end{document}